\documentclass[10pt,leqno]{article}
\usepackage{amsmath}
\usepackage{amsfonts} 
\usepackage{amsthm}
\usepackage{indentfirst} 

\theoremstyle{plain}
\newtheorem{theorem}{\indent\rm T\,h\,e\,o\,r\,e\,m\;}[section]
\newtheorem{lemma}{\indent\rm L\,e\,m\,m\,a\;}[section]

\theoremstyle{definition}

\theoremstyle{remark}

\makeatletter                                                  %
\renewcommand*{\@seccntformat}[1]{
  \csname the#1\endcsname\;-                                   %
}                                                              %
\renewcommand{\section}{\@startsection{section}{1}{0mm}        %
   {1.5\baselineskip}
   {1\baselineskip}
   {\indent\normalfont\normalsize\bfseries}
   }                                                           %
\renewcommand*{\@seccntformat}[1]{
  \normalfont\bfseries\csname the#1\endcsname\;-               %
}                                                              %
\renewcommand\subsection{\@startsection                        %
  {subsection}{2}{0mm}
  {1.5\baselineskip}
  {1\baselineskip}
  {\indent\normalfont\normalsize\itshape}}
\renewcommand*{\@seccntformat}[1]{
  \normalfont\bfseries\csname the#1\endcsname\;-               %
}                                                              %
\renewcommand\subsubsection{\@startsection                     %
  {subsubsection}{2}{0mm}
  {1.5\baselineskip}
  {1\baselineskip}
  {\indent\normalfont\normalsize\texttt}}
\makeatother                                                   %
\renewcommand{\qedsymbol}{$\square$}

\newcommand{\dx}{\mathrm{d}}
\newcommand{\eps}{\varepsilon} 
\newcommand{\R}{\mathbb{R}}
\newcommand{\N}{\mathbb{N}}

\newcommand{\C}{\mathbb{C}}

\newcommand{\Stilde}{\widetilde{S}}

\newcommand{\second}{\prime\prime}
\newcommand{\sameorder}{\asymp}
\newcommand{\Odip}[2]{\mathcal{O}_{#1}\!\left(#2\right)\mathchoice{\!}{}{}{}}

\newcommand{\Odipm}[2]{\mathcal{O}_{#1}\bigl(#2\bigr)\mathchoice{\!}{}{}{}}
\newcommand{\Odig}[1]{\mathcal{O}\Bigl(#1\Bigr)\mathchoice{\!}{}{}{}}

\newcommand{\Odim}[1]{\mathcal{O}\bigl(#1\bigr)}
\newcommand{\Odi}[1]{\Odip{}{#1}}
\newcommand{\odip}[2]{{o}_{#1}\!\left(#2\right)\mathchoice{\!}{}{}{}}
\newcommand{\odi}[1]{\odip{}{#1}}
 
\makeatletter
\newenvironment{Biggcases}{%
  \matrix@check\Biggcases\env@Biggcases
}{%
  \endarray %
}
\def\env@Biggcases{%
  \let\@ifnextchar\new@ifnextchar
  \Biggl\lbrace
  \def\arraystretch{1.2}%
  \array{@{}l@{\quad}l@{}}%
}
\makeatother

\begin{document}
\thispagestyle{empty}

\vskip -8in
\begin{center}
 \rule{8.5cm}{0.5pt}\\[-0.1cm] {\small Riv.\, Mat.\, Univ.\, Parma,\,
Vol. {\bf x} \,(2016), \,000-000}\\[-0.25cm] \rule{8.5cm}{0.5pt}
\end{center}
\vspace {2.2cm}

\begin{center}
{\sc\large Alessandro Languasco}  
\end{center}
\vspace {1.5cm}

\centerline{\large{\textbf{Applications of  some exponential sums on prime powers: a survey\footnote{Submitted for the Proceedings of the ``Terzo Incontro Italiano di Teoria dei Numeri", Scuola Normale Superiore, Pisa, 21-25 Settembre 2015.}
\footnote{This research was partially supported by the grant PRIN2010-11 \textsl{Arithmetic Algebraic Geometry and Number Theory.}}
}}}

\renewcommand{\thefootnote}{\fnsymbol{footnote}}


\renewcommand{\thefootnote}{\arabic{footnote}}
\setcounter{footnote}{1}

\vspace{1,5cm}
\begin{center}
\begin{minipage}[t]{10cm}

\small{ \noindent \textbf{Abstract.} 
Let $\Lambda$ be the von Mangoldt function
and  $N,\ell\geq 1$ be two integers.
We will see some results by the author and Alessandro Zaccagnini obtained using the original 
Hardy \& Littlewood  circle method function, {\it i.e.}
\begin{equation*} 
\widetilde{S}_{\ell}(\alpha)
=
\sum_{n=1}^{\infty} 
\Lambda(n) e^{-n^{\ell}/N}
e(n^{\ell}\alpha),
\end{equation*}
where $e(x)=\exp(2\pi i x)$, instead of 
\(
S_{\ell}(\alpha)
=
\sum_{n=1}^{N} 
\Lambda(n)  
e(n^{\ell}\alpha)
\).
We will also motivate why, for some short interval additive problems, the approach 
using $\widetilde{S}_{\ell}(\alpha)$ gives sharper results than the ones that can be obtained 
with $S_{\ell}(\alpha)$.
The final section of this paper is devoted to correct an oversight occurred in   \cite{LanguascoZ2013a}
 and  \cite{LanguascoZ2015a}. 
\medskip

\medskip

\noindent \textbf{Mathematics~Subject~Classification~(2010):} Primary 11P32; Secondary 11P55, 11P05, 44A10, 33C10.} 
\end{minipage}
\end{center}


%
%
%
%
%
%
%
%
%
%
%
%

 \section{Introduction}
 In a series of recent papers in collaboration with Alessandro Zaccagnini, we
 proved several results about additive problems with primes and prime powers.
 Our main tool was the circle method in its original form used 
 by Hardy \& Littlewood \cite{HardyL1923}, since, for these applications, it lets us get  
 stronger results than the ones obtainable by the more recent
 setting. This means to work with the function
\begin{equation}
\label{tildeS-def}
\widetilde{S}_{\ell}(\alpha)
=
\sum_{n=1}^{\infty} 
\Lambda(n) e^{-n^{\ell}/N}
e(n^{\ell}\alpha),
\end{equation}
where  $N,\ell\geq 1$ are two integers, $\Lambda$ is the von Mangoldt function and $e(x)=\exp(2\pi i x)$, instead
of 
\begin{equation}
\label{S-def}
S_{\ell}(\alpha)
=
\sum_{n=1}^{N} 
\Lambda(n)  
e(n^{\ell}\alpha).
\end{equation}
 
 We see now the statements we were able to prove
 and we compare them with the known results obtained
 using \eqref{S-def}. Then we will try to give a motivation
 why, for such kind of problems, using \eqref{tildeS-def}
 is better than using \eqref{S-def}.

\section{Statements} 
\label{statements}

All the following results are  asymptotic formulae
concerning suitable averages 
for functions counting the number of ways an
integer is representable as a sum of primes, prime powers
or powers. It is well known that usually the behaviour 
of such functions is erratic but averaging them
is, in many cases, sufficient to get much more regular
problems.  
We also remark that every asymptotic formula we will write below has,
as a consequence, the existence of an integer
representable in the way stated in the result in  intervals
of the mentioned size.
 
Our main goal is to get, both assuming the Riemann Hypothesis
and unconditionally, an asymptotic relation involving zeros of the Riemann
zeta-function $\zeta(s)$, which either holds in the shorter known interval
or has the best known error estimate.

We start with the most famous additive problem: the Goldbach one.
Letting
 \[
 R_{G}(n)= \sum_{m_1 + m_2  = n} \Lambda(m_1)  \Lambda(m_2),
 \]  
 we obtained
\begin{theorem}[Languasco-Zaccagnini \cite{LanguascoZ2012a}]
\label{Bhowmik-thm}
Let $N$ be a sufficiently large integer.
Assuming the  Riemann Hypothesis (RH) holds, we have
\begin{equation}
\label{bhowmik-improved-estim}
\sum_{n=1}^{N} 
R_G(n) 
=
\frac{N^{2}}{2}
-2
\sum_{\rho} \frac{N^{\rho + 1}}{\rho (\rho + 1)}
+
\Odim{N \log^{3}N
},
\end{equation}
where $\rho=1/2+i\gamma$ runs over
the non-trivial zeros of  $\zeta(s)$. 
\end{theorem}
Using  \eqref{S-def},
Bhowmik \& Schlage-Puchta \cite{BhowmikS2010} proved that the error term
is $\Odim{N \log^{5}N}$. Quite recently Goldston \& Yang \cite[Theorem 1]{GoldstonY2016}, inserting
a suitable average in Bhowmik \& Schlage-Puchta's approach (see Lemma 5 of \cite{GoldstonY2016}), 
were able to give another proof of Theorem \ref{Bhowmik-thm}.
In \cite{LanguascoZ2012a} we also proved a short interval version of
Theorem \ref{Bhowmik-thm} and its analogue with  the case of representing an integer as a sum
of $k\geq 3$ primes. 
 
Inserting the Ces\`aro  (or Riesz)  weight in the picture,
a tool widely used in Number Theory history, see, \emph{e.g.}, \cite{ChandrasekharanN1961},    
we were also able to push forward the investigation
about terms of secondary order in Theorem \ref{Bhowmik-thm}  thus getting  
\begin{theorem}[Languasco-Zaccagnini  \cite{LanguascoZ2015a}]
\label{CesaroG-thm}
Let $N$ be a sufficiently large integer and
 $k>1$ be a real number. Then 
\begin{align}
\notag
  \sum_{n \le N} R_G(n) \frac{(1 - n/N)^k}{\Gamma(k + 1)}
  &=
  \frac{N^{2}}{\Gamma(k + 3)}
  -
  2
  \sum_{\rho} \frac{\Gamma(\rho)}{\Gamma(\rho + k + 2)} N^{\rho+1} \\
  &\qquad+
\label{expl-form-Goldbach}
  \sum_{\rho_1} \sum_{\rho_2}
    \frac{\Gamma(\rho_1) \Gamma(\rho_2)}{\Gamma(\rho_1 + \rho_2 + k + 1)}
    N^{\rho_1 + \rho_2}
  +
  \Odipm{k}{N} ,
\end{align}
where $\rho$, with or without subscripts, runs over
the non-trivial zeros of  $\zeta(s)$
and $\Gamma$ is the Euler Gamma-function.
\end{theorem}
In the original version of Theorem \ref{CesaroG-thm} the error term was erroneously 
written as  $\Odipm{k}{N^{1/2}}$; see Section \ref{corrigendum} below for the correction.

 Goldston \& Yang \cite[Theorem 2]{GoldstonY2016}, assuming RH, were recently able to prove that
 \begin{equation}
 \label{Goldston-Yang}
   \sum_{n \le N}  R_G(n) \Bigl(1 - \frac{n}{N}\Bigr)
  =
  \frac{N^{2}}{6}
  -
  2
  \sum_{\rho} \frac{N^{\rho+1}}{\rho(\rho+1)(\rho+2)}   +  \Odim{N}
 \end{equation}
which would correspond to what \eqref{expl-form-Goldbach} 
would imply if one could take  $k=1$ in Theorem \ref{CesaroG-thm}. 

A  result similar to Theorem \ref{CesaroG-thm} can be obtained for the problem of
representing an integer as a sum of a prime and a square, {\it i.e.}, the so-called
Hardy-Littlewood numbers since in \cite{HardyL1923} they conjectured
that every large non-square integer   should be representable
in such a way. Letting
\[
  R_{\textit{HL}}(n)  =  \sum_{m_1 + m_2^2 = n} \Lambda(m_1),
\]
we proved
 \begin{theorem}[Languasco-Zaccagnini  \cite{LanguascoZ2013a}]
 \label{CesaroHL-thm}
 Let $N$ be a sufficiently large integer and
  $k>1$ be a real number. Then 
 
\begin{align*}
  \sum_{n \le N} R_{\textit{HL}}(n)  \frac{(1 - n/N)^k}{\Gamma(k + 1)}
  &=
  \frac{\pi^{1 / 2}}2 \frac{N^{3 / 2}}{\Gamma(k + 5 / 2)}
  -
  \frac 12 \frac{N}{\Gamma(k + 2)}
  \\
  &
  -
  \frac{\pi^{1 / 2}}2
  \sum_{\rho}
    \frac{\Gamma(\rho)}{\Gamma(k + 3 / 2 + \rho)} N^{1 / 2 + \rho} 
     \end{align*}
 \begin{align*}
 \hskip2cm 
 &+
  \frac 12
  \sum_{\rho}
    \frac{\Gamma(\rho)}{\Gamma(k + 1 + \rho)} N^{\rho}
  +
  \frac{N^{3 / 4 - k / 2}}{\pi^{k + 1}}
  \sum_{\ell \ge 1}
    \frac{J_{k + 3 / 2} (2 \pi \ell N^{1 / 2})}{\ell^{k + 3 / 2}} \\
  &-
  \frac{N^{1 / 4 - k / 2}}{\pi^k}
  \sum_{\rho} \Gamma(\rho) \frac{N^{\rho / 2}}{\pi^\rho}
    \sum_{\ell \ge 1}
      \frac{J_{k + 1 / 2 + \rho} (2 \pi \ell N^{1 / 2})}
           {\ell^{k + 1 / 2 + \rho}}
  +
  \Odipm{k}{N^{1/2}},
\end{align*}
where   $\rho$  runs over
the non-trivial zeros of $\zeta(s)$, $\Gamma$ is the Euler Gamma-function
and $J_{\nu} (u)$ denotes the Bessel function of complex order $\nu$
and real argument $u$.
\end{theorem}
In the original version of Theorem \ref{CesaroHL-thm} the error term was erroneously 
written as  $\Odipm{k}{1}$; see Section \ref{corrigendum} below for the correction. 
We recall here  the
Sonine representation of such a Bessel function:
\begin{equation}
\label{Bessel-def}
  J_\nu(u)
  :=
  \frac{(u / 2)^\nu}{2 \pi i}
  \int_{(a)} s^{- \nu - 1} e^s e^{- u^2 / 4 s} \, \dx s,
\end{equation}
where $a > 0$ and $u,\nu \in \C$ with $\Re(\nu) > -1$.

Theorem \ref{CesaroHL-thm}  has no analogue proved using  \eqref{S-def}.
Moreover it seems that it was the first time a complex order 
Bessel function is involved in additive results. In particular
these terms come from the modularity of the complex $\theta$-function
$\theta(z)  =   \sum_{m = -\infty}^{+\infty} e^{-m^2 z}$, $z\in \C$, given by the formula
\begin{equation}
\label{theta-functeq}
  \theta(z)
  =
  \Bigl( \frac \pi z \Bigr)^{1/2}
  \theta\Bigl( \frac{\pi^2} z \Bigr)\quad
  \textrm{for} \ \Re (z)>0.
\end{equation}

We then focused our attention on problems involving prime squares.
In 1938, Hua \cite{Hua1938} considered the problem of representing an integer
as a sum of a prime and of two prime squares. Letting
${\mathcal A} = \{ n\in \N : n\equiv 1 \bmod 2;\ n \not \equiv 2\bmod 3\}$,
Hua proved that the set of integers  $n\in [1,N]\cap  {\mathcal A}$
that cannot be representable in such a way has a cardinality $\ll N/\log ^A N$
for some $A>0$. Recently  L.~Zhao  \cite{Zhao2014a} improved this
estimate to $\ll N^{1/3+\eps}$, for every  $\eps>0$. 
Letting
\[
  r(n)
  =
  \sum_{p_{1} + p_{2}^2 + p_{3}^{2} = n} \log p_{1} \log p_{2} \log p_{3},
\]
we proved the following
\begin{theorem}[Languasco-Zaccagnini \cite{LanguascoZ2016a}]
\label{Hua-thm} 
Assume RH. Then
\[
   \sum_{n = N+1}^{N + H}   r(n) = \frac{\pi}{4} HN
   +
\Odig{H^{1/2}N (\log N)^{2} +H N^{3/4}(\log N)^{3}   +   H^{2}(\log N)^{3/2}} 
 \] 
as $N \to \infty$,
uniformly for  $\infty((\log N)^{4})\le H \le \odi{N(\log N)^{-3/2}}$, where $f=\infty(g)$ means $g=\odi{f}$. 
\end{theorem}

We were also able to prove an unconditional version of Theorem \ref{Hua-thm} 
that holds just for  $H \ge N^{7/12+\eps}$.
Using \eqref{S-def}, a similar result to Theorem \ref{Hua-thm} 
holds just for $H \ge N^{1/2+\eps}$ under RH; this is due to a limitation I'll discuss 
in Section \ref{motivations} below, see the remark after Lemma \ref{App-BCP-Gallagher}. 


 Similar, but harder, problems with prime and prime squares are the binary ones.
 Considering
\[ 
r^{\second}_{1,2}(n) = 
 \sum_{p_{1}+p_{2}^2=n} 
\log p_{1} \log p_{2}\, ,
\]
we got the following
\begin{theorem}[Languasco-Zaccagnini \cite{LanguascoZ2015c}]
\label{density-32}
Assume RH. Then
\[
\sum_{n = N+1}^{N + H} r^{\second}_{1,2}(n)
=
H N^{1/2} 
+ 
 \Odig{\frac{H^2}{N^{1/2}}+  N^{3/4}(\log N)^3+ HN^{1/3}(\log N)^2}
\]
as $N\to \infty$ uniformly for $\infty(N^{1/4}(\log N)^{3}) \le H \le \odi{N}$. 
\end{theorem}

We also obtained
 that the error term $\Odig{H \exp \Big( - C  \Big( \frac{\log N}{\log \log N} \Big)^{1/3} \Big)}$
unconditionally holds uniformly for $N^{7/12+\eps}\le H \le N^{1-\eps}$.
Recently, Suzuki \cite{Suzuki2015} proved that this result holds in $N^{1/2}\exp \Big( - C  \Big( \frac{\log N}{\log \log N} \Big)^{1/3} \Big)\le H \le N^{1-\eps}$.

As far as we know, Theorem \ref{density-32} is the first short interval result for this problem;
we recall that the method in Plaksin \cite{Plaksin1981} (Lemma 11 there) gives long interval results only
and that an adaptation his idea to the short interval $[N,N+H]$ leads to   weaker results.
We also remark that using   \eqref{S-def} we can only reach the level 
$H \ge N^{1/2+\eps}$.
Finally we   recall that in \cite{LanguascoZ2015c} we also have results 
for the $p+m^{2}=n$ case.

Even harder  problems are the ones involving two prime squares, or a prime square and a square.
As far as we know, for these problems just long interval averages were available, see 
the Introduction of the paper by Daniel \cite{Daniel2001}. 
Letting
 \[
r^{\second}_{2,2}(n) = 
\sum_{p_{1}^{2}+p_{2}^2=n}
\log p_1 \log p_2
\quad
\textrm{and}
\quad
r^{\prime}_{2,2}(n) = 
\sum_{p^{2}+m^2=n}
\log p \, ,
\]
we obtained the following
\begin{theorem}[Languasco-Zaccagnini \cite{LanguascoZ2015d}]
\label{two-prime-squares}
Assume RH. Then
\[
\sum_{n = N+1}^{N + H} r^{\second}_{2,2}(n)
=
 \frac{\pi}{4}   H
+ \Odig{\frac{H^2}{N} +  H^{1/2}N^{1/4}(\log N)^{3/2}}
\]
as  $N\to \infty$ uniformly for $\infty(N^{1/2}(\log N)^{3})\le H \le  \odi{N}$. 
Moreover we also have that
\[
\sum_{n = N+1}^{N + H} r^{\prime}_{2,2}(n)
=
\frac{\pi}{4} H 
+
 \Odig{\frac{H^2}{N}+\frac{H \log \log N} {(\log N)^{1/2}}  }
 \]
as  $N\to \infty$ uniformly for $\infty(N^{1/2}(\log N)^{2}) \le H \le \odi{N}$.  
\end{theorem}
We also obtained
 that the error term $\Odig{H \exp \Big( - C  \Big( \frac{\log N}{\log \log N} \Big)^{1/3} \Big)}$
unconditionally holds uniformly for $N^{7/12+\eps}\le H \le N^{1-\eps}$ in both cases.

\section{Settings for Theorems \ref{Bhowmik-thm} and \ref{Hua-thm}-\ref{two-prime-squares}}
\label{setting}

The proofs of Theorems \ref{Bhowmik-thm} and \ref{Hua-thm}-\ref{two-prime-squares}  
all work with similar settings
but, for simplicity, here we just write the Goldbach case in short interval. 

Let  $N,H>0$,  $H=\odi{N}$, be  integers, 
\begin{equation}
\label{z-def}
z= 1/N-2\pi i\alpha
\end{equation} 
and $M_1(\alpha) = 1/z$ be
the expected main term for the involved exponential sum, see Lemma  \ref{Linnik-lemma} below. Moreover let
\[
U(\alpha,H) = \sum_{1\leq m\leq H}e(m\alpha).
\]
We have 
\begin{align}
\notag
\sum_{n=N+1}^{N+H}
e^{-n/N}R_G(n) 
&=
\int_{-1/2}^{1/2} M_1(\alpha)^2 U(\alpha,H)\, \dx \alpha 
\\
\label{main-circle-setting}
&
+
\int_{-1/2}^{1/2} (\widetilde{S}_1(\alpha)^2 -M_1(\alpha)^2)U(\alpha,H)\, \dx \alpha.
\end{align}
The first term on the right hand side will lead to the expected main term
while the second one is connected with a $L^2$-estimate for 
$\widetilde{S}_1(\alpha)-M_1(\alpha)$.

In fact, using the algebraic relation $f^2-g^2=(f-g)^2+2g(f-g)$, the Cauchy-Schwarz inequality and 
$
U(\alpha,H)\ll  \min (H; \vert \alpha\vert ^{-1}),
$
the main quantity to get 
an upper bound for the second term on the right hand 
side of \eqref{main-circle-setting}  is  
\begin{equation}
\label{split}
H
\int_{-1/H}^{1/H} \vert\widetilde{S}_1(\alpha) -M_1(\alpha)\vert^2 \, \dx \alpha
+
\int_{1/H}^{1/2} \vert\widetilde{S}_1(\alpha) -M_1(\alpha)\vert^2\,  \frac{\dx \alpha}{\alpha},
\end{equation}
since a mean-square estimate of $\widetilde{S}_1(\alpha)$,
needed to estimate the mixed term in the previously mentioned  
agebraic relation, can be obtained
using the Prime Number Theorem.
A partial integration reveals that the key quantity is then
\[
\int_{-\xi}^{\xi} \vert\widetilde{S}_1(\alpha) -M_1(\alpha)\vert^2 \, \dx \alpha
\]
with $0<\xi\leq 1/2$.
Hence we  clearly  need results on computing the main term
and on such a truncated mean-square. For this second object, it
is clear that the larger $\xi$-uniformity  we have, the stronger result
on $H$ we get.

In the next section we will see that, for $\ell \ge 2$, the  $\xi$-uniformity for the $\Stilde_{\ell}(\alpha)$
 truncated $L^2$-average is, at least in the conditional case,  much larger than the corresponding one for 
 $S_{\ell}(\alpha)$.

We finally remark that in many cases it is easy to remove the contribution 
of the exponential weight  $e^{-n/N}$ from \eqref{main-circle-setting} since 
$e^{-n/N}= e^{-1} +\Odi{H/N}$ for $n\in [N+1,N+H]$. For proving the result stated in 
Theorem \ref{Bhowmik-thm}
such a  relation for $e^{-n/N}$ was not good enough and we had to modify the 
previously described setting to be able to detect the secondary 
order terms in \eqref{bhowmik-improved-estim},
see eq. (3.1) of the proof of Theorem 3.1 in \cite{LanguascoZ2012a}.
 
\section{General results on $\Stilde_{\ell}(\alpha)$ and $S_{\ell}(\alpha)$ and motivations} 
\label{motivations}

For \eqref{tildeS-def} we have an explicit formula 
which can be proved via a Mellin transform argument. 
It generalizes and slightly sharpens what Linnik \cite{Linnik1946}-\cite{Linnik1952} proved for $\ell =1$. Recall that $z= 1/N-2\pi i\alpha$.

\begin{lemma}[Lemma 2 of  \cite{LanguascoZ2016a}]
\label{Linnik-lemma}
Let $\ell \ge 1$ be an integer, $N \ge 2$  and $\alpha\in [-1/2,1/2]$.
Then
\begin{equation}
\Stilde_{\ell}(\alpha)  
= 
\frac{\Gamma(1/\ell)}{\ell z^{1/\ell}}
- 
\frac{1}{\ell}\sum_{\rho}z^{-\rho/\ell}\Gamma\Bigl(\frac{\rho}{\ell}\Bigr) 
+
\Odip{\ell}{1},
\end{equation}
where $\rho=\beta+i\gamma$ runs over
the non-trivial zeros of $\zeta(s)$.
\end{lemma} 

This is in fact an explicit formula for $\Stilde_{\ell}(\alpha)$ and plays a key role
in every result listed in  Section \ref{statements}; unfortunately
a direct analogue of such a result for \eqref{S-def} is unknown. 
Lemma \ref{Linnik-lemma} means that the expected main term comes form the 
$\frac{\Gamma(1/\ell)}{\ell z^{1/\ell}}$ part. In fact  it can be 
usually obtained via the following 
\begin{lemma} [Lemma 4 of \cite{LanguascoZ2016a}]
 \label{Laplace-formula}
Let $N$ be a positive integer and 
$\mu > 0$.
Then
\[
  \int_{-1 / 2}^{1 / 2} z^{-\mu} e(-n \alpha) \, \dx \alpha
  =
  e^{- n / N} \frac{n^{\mu - 1}}{\Gamma(\mu)}
  +
  \Odip{\mu}{\frac{1}{n}},
\]
uniformly for $n \ge 1$.
\end{lemma}
The proof of Lemma \ref{Laplace-formula} directly follows from the 
Laplace formula \eqref{Laplace-transf}.

As we have seen in Section \ref{setting} we now need information
on some truncated  $L^{2}$ average.   
For \eqref{tildeS-def} we have  the following  result   
 (the case $\ell=1$ under RH is in   \cite{LanguascoP1994}).
\begin{lemma}[Lemma 3 of \cite{LanguascoZ2016a} and 
Lemma 1 of \cite{LanguascoZ2015d}]
 \label{LP-Lemma-gen} 
Let $\eps$ be an arbitrarily small
positive constant,  $\ell \ge 1$ be an integer, $N$ be a
sufficiently large integer and $L=\log N$. Then there exists a positive constant 
$c_1 = c_{1}(\eps)$, which does not depend on $\ell$, such that 
\[
\int_{-\xi}^{\xi} \,
\Bigl\vert
\Stilde_\ell(\alpha) - \frac{\Gamma(1/\ell)}{\ell z^{1/\ell}}
\Bigr\vert^{2}
\dx \alpha 
\ll_{\ell}
 N^{2/\ell-1} \exp \Big( - c_{1}  \Big( \frac{L}{\log L} \Big)^{1/3} \Big) 
\]
uniformly for $ 0\le \xi < N^{-1 +5/(6\ell) - \eps}$.
Assuming RH we get 
\[
\int_{-\xi}^{\xi} \,
\Bigl\vert
\Stilde_\ell(\alpha) - \frac{\Gamma(1/\ell)}{\ell z^{1/\ell}}
\Bigr\vert^{2}
\dx \alpha 
\ll_{\ell}
N^{1/\ell}\xi L^{2}
\]
uniformly  for  $0 \le \xi \le 1/2$.
\end{lemma} 
The proof of Lemma \ref{LP-Lemma-gen}, even if it is a direct one, is quite lengthy
due to some tedious computation to be performed for handling the contribution
of   complex powers of the  $z^\rho$-type.

In Theorems 3.1-3.2 of   \cite{LanguascoZ2013b}  (see also Lemma 1 of   \cite{LanguascoZ2016a}
and  the Theorem in Languasco \cite{Languasco1998b})
we proved a direct analogue of the previous lemma for \eqref{S-def}.
Letting
\[
 T_{\ell}(\alpha)  
 =
 \sum_{n=1}^{N}  
e(n^{\ell}\alpha)\, ,
\]
which now plays the role of the expected main term for $ S_{\ell}(\alpha)$,
we have
\begin{lemma}  
\label{App-BCP-Gallagher}
Let $\eps$ be an arbitrarily small
positive constant, $\ell  > 0$ be a real number,  $N$ be a
sufficiently large integer and $L=\log N$. Then there exists a positive constant 
$c_1 = c_{1}(\eps)$, which does not depend on $\ell$, such that
\[
\int_{-\xi}^{\xi}
\vert
S_{\ell}(\alpha) - T_{\ell}(\alpha)  
\vert^2 
\, \dx \alpha
\ll_{\ell}   N^{2/\ell -1}
\Bigl(
\exp \Big( - c_{1}  \Big( \frac{L}{\log L} \Big)^{1/3} \Big)
+
\frac{L^{2}}{\xi N}
\Bigr),
\]
uniformly for  $1/(2N) \le \xi  < N^{-1 +5/(6\ell) - \eps}$.  
Assuming   RH we get 
\[
\int_{-\xi}^{\xi}
\vert
S_{\ell}(\alpha) - T_{\ell}(\alpha) 
\vert^2 
\, \dx \alpha
\ll_{\ell} 
 N^{1/\ell} \xi L^{2}  + \frac{N^{2/\ell-2} L^{2}}{\xi},
\]
uniformly for  $1/(2N) \le \xi \le N^{1/\ell-1}/2$.
For $\ell=1$ the last estimate can be replaced by 
$N\xi L^{2} + \min\Bigl(\frac{L^{2}}{\xi(\log 3\xi)^2};N\xi L^4\Bigr)$
uniformly for  $1/(2N) \le \xi \le 1/2$.
\end{lemma}

We remark that the $\xi$-uniformity here is, at least in the conditional case and for $\ell\ge 2$,
worse than the one obtained in Lemma \ref{LP-Lemma-gen}. This
is in fact a key limitation for the kind of problems we are working on.
We try now to explain why we  encountered it.
The proof  of Lemma \ref{App-BCP-Gallagher}  is obtained via Gallagher's lemma \cite{Gallagher1970}
which gives
\begin{align}
\notag
 \int_{-\xi}^{\xi}
&\vert
S_{\ell}(\alpha) - T_{\ell}(\alpha) 
\vert^2 
\, \dx \alpha
=
\int_{-\xi}^{\xi}
\Big \vert
 \sum_{N \leq n^{\ell} \leq  2N} ( \Lambda(n)-1 )e(n^{\ell} \alpha)
\Big\vert^{2}
\dx \alpha
\\
\label{Gallagher}
&
\ll  
\xi^{2}
\int_{N}^{2N-H}
\Bigl(
\sum_{\substack{x \leq n^{\ell} \leq x + H \\ N \leq n^{\ell} \leq  2N}}
 ( \Lambda(n)-1)
\Bigr)^{2}
\dx x  + \Odim{H^{3}N^{2/\ell-2} L^{2}
+ H}\, ,
\end{align}
where we defined 
$H=1/(2\xi)$.
From $x \leq n^{\ell} \leq x + H$ one immediately gets
$x^{1/\ell} \leq n \leq (x + H)^{1/\ell}$ and hence, to have
$(x + H)^{1/\ell} -x^{1/\ell}\gg 1$,  we need
$ H x^{1/\ell-1} \gg 1$, that is, essentially,
$ \xi \ll N^{1/\ell-1}$.
 
The estimate of the integral at the right hand side of \eqref{Gallagher}
can be performed exploiting density estimates for the non-trivial
zeros of the Riemann-zeta function, as in Saffari \& Vaughan \cite{SaffariV1977}.
We further notice that if  $\ell \geq 2$ and $\xi \ge N^{1 / \ell-1}$, then the bound $\ll \xi^2 N \log N$ 
(which in this range is worse than $\xi N^{1/\ell} \log^2 N$)
follows immediately from the Prime Number Theorem.

For $\ell=1$  and assuming RH,  a slightly more careful reasoning
 (see Languasco \cite{Languasco1998b}) leads to the alternative estimate in  Lemma 
 \ref{App-BCP-Gallagher}. This
reveals that  Lemmas \ref{LP-Lemma-gen}  and \ref{App-BCP-Gallagher} in this case
have the same the $\xi$-uniformity and they
differ  just for some logarithmic factors; in fact the gain in our Theorem \ref{Bhowmik-thm}
directly follows from such a slight improvement.\footnote{Goldston \& Yang \cite{GoldstonY2016}, introducing a further external average technique (see their Lemmas 5 and 6), were able to overcome this problem at least for the $\ell=1$ case.}

For $\ell \geq 2$  and assuming RH,   Lemmas \ref{LP-Lemma-gen}  and \ref{App-BCP-Gallagher}
have quite different uniformities in the $\xi$-aspect and 
this means that additive problems with prime powers summands fit better with the \eqref{tildeS-def}-approach
since the larger $\xi$ we get, the shorter $H$ we can reach in \eqref{split} and hence in \eqref{main-circle-setting}.

However,  Lemma \ref{App-BCP-Gallagher} is quite useful in proving Diophantine
approximation results using linear forms in prime and prime powers
see, \emph{e.g.},  \cite{LanguascoS2012},
 \cite{LanguascoZ2010b}, \cite{LanguascoZ2012c},
\cite{LanguascoZ2012e},
\cite{LanguascoZ2013b}.

For problems involving prime squares and squares we also needed information
 about some truncated $L^2$ estimates for $\Stilde_\ell(\alpha)$ and 
 about
\begin{equation}
\label{omega-def}
\omega_{\ell}(\alpha)
=
\sum_{m=1}^{\infty}   
e^{-m^{\ell}/N} e(m^{\ell}\alpha)
=
\sum_{m=1}^{\infty}   
e^{-m^{\ell}z} .
\end{equation}
Using a well-known Montgomery \& Vaughan result \cite{MontgomeryV1974}, we proved  
\begin{lemma}[Lemma 2 of \cite{LanguascoZ2015d}]
%
Let  $0<\xi\leq 1/2$, $\ell\ge 2 $ be an integer,
$N$ be a
sufficiently large integer and $L=\log N$. Then
\[
\int_{-\xi}^{\xi} 
|\omega_{\ell}(\alpha)|^2 \ \dx\alpha 
\ll_{\ell}
  \xi N^{1/\ell}
  + 
\begin{Biggcases}
L & \text{if}\ \ell =2\\
1 & \text{if}\ \ell > 2
\end{Biggcases}\]
and
\[
\int_{-\xi}^{\xi} 
|\Stilde_{\ell}(\alpha)|^2 \ \dx\alpha 
\ll_{\ell}
\xi N^{1/\ell} L  +
\begin{Biggcases}
L^{2} & \text{if}\ \ell =2\\
1 & \text{if}\ \ell > 2.
\end{Biggcases}
\]
\end{lemma}
Moreover in some cases we also needed information 
about the fourth moment for $\Stilde_2(\alpha)$. We proved the following
\begin{lemma} [Lemma 5 of \cite{LanguascoZ2016a}]
\label{Rieger}
Let $N$ be a
sufficiently large integer and $L=\log N$. We have 
\[
\int_{-1/2}^{1/2} \,
\vert
\Stilde_2(\alpha) 
\vert^{4}
\, \dx \alpha 
\ll
N L^{2}.
\]
\end{lemma} 
The proof of Lemma \ref{Rieger} exploits cancellation and uses a  Rieger's result
(see Satz 3 on page 94 of
\cite{Rieger1968}) about the number of solutions of the $p_1^2+p_2^2 = p_3^2+p_4^2$
problem.
 
We conclude this section remarking another limitation we encountered
 using \eqref{S-def}. In problems involving 
two prime squares we have to approximate  
$T_{2}(\alpha) = \sum_{1 \le m^{2} \le  N} e(m^{2} \alpha )$
with $f_{2}(\alpha) =(1/2)\sum_{1 \leq m\leq N} m^{-1/2}e(m\alpha)$;
a standard approach in Waring-Goldbach problems.
The main term for the $p_1^2+p_2^2=n$ problem has order of magnitude 
$\asymp H$; this comes from the evaluation of the integral of 
$f_{2}(\alpha)^{2}U(-\alpha,H)$ which 
can be computed using  Lemma 2.9 of Vaughan \cite{Vaughan1997}.
But we also need to evaluate the quantity
\[
\Bigl\vert   
\int_{-1/H}^{1/H}  (T_{2}(\alpha)^{2}-f_{2}(\alpha)^{2})U(-\alpha,H)
 e(-N\alpha)\, \dx \alpha
 \Bigr\vert.
\]
By 
Theorem 4.1 of Vaughan \cite{Vaughan1997} we know that
$\vert T_{2}(\alpha) - f_{2}(\alpha) \vert \ll (1+\vert \alpha \vert N)^{1/2}$
and hence it is not hard to see that the previous integral is
$\ll NH^{-1/2}$.
Since the expected order of magnitude of the 
main term is $H$, the previous estimate
is under control  if and only if $H\ge N^{2/3+\eps}$;
a result which is weaker 
than the one we obtained, see the statement of Theorem \ref{two-prime-squares}
above.

 \section{Settings for Theorems \ref{CesaroG-thm}-\ref{CesaroHL-thm}}
 
 The method in this case is based on  Laplace transforms. In particular we used  a formula due to Laplace
\cite{Laplace1812}, namely
\begin{equation}
\label{Laplace-transf}
  \frac 1{2 \pi i}
  \int_{(a)} v^{-s} e^v \, \dx v
  =
  \frac1{\Gamma(s)},
\end{equation}
where $\Re(s) > 0$ and $a > 0$, see, e.g., formula 5.4(1) on page 238
of \cite{ErdelyiMOT1954a}.
In the following we will need the general case of \eqref{Laplace-transf}
which, {\it e.g.}, can be found in   \cite{Azevedo2002}, 
formulae (8)-(9):
\begin{equation}
\label{Laplace-eq-1}
  \frac1{2 \pi}
  \int_{\R} \frac{e^{i D u}}{(a + i u)^s} \, \dx u
  =
  \begin{cases}
    \dfrac{D^{s - 1} e^{- a D}}{\Gamma(s)}
    & \text{if $D > 0$,} \\
    0
    & \text{if $D < 0$,}
  \end{cases}
\end{equation}
which is valid for $\sigma = \Re(s) > 0$ and $a \in \C$ with
$\Re(a) > 0$, and
\begin{equation}
\label{Laplace-eq-2}
  \frac1{2 \pi}
  \int_{\R} \frac 1{(a + i u)^s} \, \dx u
  =
  \begin{cases}
    0     & \text{if $\Re(s) > 1$,} \\
    1 / 2 & \text{if $s = 1$,}
  \end{cases}
\end{equation}
for $a \in \C$ with $\Re(a) > 0$ (for the case $s=1$ the principal value of the integral has to be taken).
Formulae \eqref{Laplace-eq-1}-\eqref{Laplace-eq-2} enable us to write
averages of arithmetical functions by means of line integrals as we
will see  below.
We recall that Walfisz, see \cite[Ch.~X]{Walfisz1957}
used \eqref{Laplace-transf} in the particular case  where $s \in \N$, $s\geq 2$.
 
As we did in Section \ref{setting}, we  just write the setting about the Goldbach problem
since the other one is similar. We need to slightly change the notation of
$\Stilde_{\ell}(\alpha)$ since now we need to let  $\Im(z)$
run over $\R$. Hence from now on we use
\[
\Stilde_{\ell}(z)
 = \sum_{n=1}^{\infty} \Lambda(n)e^{-n^{\ell}z} ,
 \quad
 z= 1/N +i y,
 \ y\in \R.
 \]

We first recall that the Prime Number Theorem  is equivalent, via
Lemma~\ref{Linnik-lemma2} below, to the statement
\begin{equation}
\label{PNT-equiv}
  \Stilde_{1}(z)
  \sim
  \frac{1}{z}
  \qquad\text{for $N \to \infty$,}
\end{equation}
for the proof see Lemma~9 in Hardy \&
Littlewood \cite{HardyL1923}.
By \eqref{tildeS-def} and \eqref{z-def} we have
\(
  \Stilde_1(z)^2
  =
  \sum_{n \ge 1} R_G(n) e^{-n z}.
\)
Hence, for $a=1/N$, $N \in \N$ with $N > 0$,  we have
\begin{equation}
\label{starting-point}
 \frac 1{2 \pi i}
  \int_{(1/N)} e^{N z} z^{- k - 1}  \Stilde_1(z)^2 \, \dx z
  =
  \frac 1{2 \pi i}
  \int_{(1/N)} e^{N z} z^{- k - 1}
    \sum_{n \ge 1} R_{G}(n) e^{- n z} \, \dx z.
\end{equation}
Since
\(
  \sum_{n \ge 1}   R_G(n) e^{- n/N}   
  =
  \Stilde_1(1/N)^2
  \sameorder
  N^{2}
\)
by \eqref{PNT-equiv}, where $f\sameorder g$ means $g \ll f \ll g$, we
can exchange the series and the line integral in
\eqref{starting-point} provided that $k>0$.
Using \eqref{Laplace-eq-1} for $n \ne N$ and \eqref{Laplace-eq-2} for
$n = N$, we see that for $k > 0$ the right hand side of \eqref{starting-point} is
\begin{align*}
  &=
  \sum_{n \ge 1} R_{G}(n)
    \Bigl[
      \frac 1{2 \pi i}
      \int_{(1/N)} e^{(N - n) z} z^{- k - 1} \, \dx z
    \Bigr]  
  =
  \sum_{n \le N} R_{G}(n) \frac{(N - n)^k}{\Gamma(k + 1)}.
\end{align*}

Summing up, for  $k > 0$ we have
\begin{equation} 
\label{fundamental}
  \sum_{n \le N}
    R_{G}(n) \frac{(N - n)^k}{\Gamma(k + 1)}
  =
 \frac 1{2 \pi i}
  \int_{(1/N)} e^{N z} z^{- k - 1} \Stilde_1(z)^2 \, \dx z,
\end{equation}
where $N \in \N$ with $N > 0$. 
This is the fundamental equation of the method since it lets us  insert the Ces\`aro (or Riesz) weight
in the picture.

Clearly  we need to work on a vertical line
and hence some of the results of Section \ref{motivations}
need to be adapted to this case. In particular we have
 \begin{lemma} [Lemma  4.1 of \cite{LanguascoZ2015a}]
\label{Linnik-lemma2}
Let $z = 1/N + iy$, where   $y \in \R$.
Then
\begin{equation*} 
  \widetilde{S}_1(z)
  =
  \frac{1}{z}
  -
  \sum_{\rho}z^{-\rho} \Gamma(\rho)
  +
  E(N,y)
\end{equation*}
where $\rho = \beta + i\gamma$ runs over the non-trivial zeros of
$\zeta(s)$ and
\begin{equation*}
  E(N,y)
  \ll
  1+
  \vert z \vert^{1/2}
  \begin{cases}
    1 & \text{if $\vert y \vert \leq 1/N$} \\
    1 +\log^2 (N\vert y\vert) & \text{if $\vert y \vert > 1/N$.}
  \end{cases}
\end{equation*}
\end{lemma}

In the original statement of Lemma \ref{Linnik-lemma2} the constant term in the
estimate for $E(N,y)$ is missing
(but the original proof is correct).  
This oversight affects the size of
the  the error term estimate for Theorems \ref{CesaroG-thm}
e \ref{CesaroHL-thm};  see Section \ref{corrigendum} below for the correction.

The remaining part of the proof of Theorem \ref{CesaroG-thm}
follows by inserting Lemma \ref{Linnik-lemma2}
into the right hand side of \eqref{fundamental}, by 
verifying that we can safely interchange sums and integrals
(this  leads to the condition $k>1$)
and  using again \eqref{Laplace-eq-1}-\eqref{Laplace-eq-2}.

For the $p+m^2$ problem, in \eqref{fundamental} one copy of $\widetilde{S}_1(z)$
has to be replaced by the function
\(
\omega_2(z)
=
\sum_{m=1}^{\infty}   
e^{-m^{2}z} .
\) 
The use of the modular relation  for $ \omega_2(z)$, which immediately
follows from \eqref{theta-functeq}, leads to 
much more complicated computations than in the Goldbach case;
in fact  integrals of the type in \eqref{Bessel-def} appeared
in several series   also
involving non-trivial zeros of the Riemann zeta-function.
Eventually, after some work, we were able to prove
that even in this case all the interchanging between series
and integrals are justified if $k>1$; thus proving Theorem \ref{CesaroHL-thm}.

In both Theorems \ref{CesaroG-thm}-\ref{CesaroHL-thm}  would be very interesting 
being able to improve the range of $k$. The most important result would be
getting  $k > 0$ but, unfortunately, this seems not to be reachable by current technologies
even if there are indications that some result of this kind should hold. For example,
the term 
\[
  \sum_{\rho_1} \sum_{\rho_2}
    \frac{\Gamma(\rho_1) \Gamma(\rho_2)}{\Gamma(\rho_1 + \rho_2 + k + 1)}
    N^{\rho_1 + \rho_2}
\]
involved in the statement of Theorem \ref{CesaroG-thm}, is in fact
absolutely convergent for $k>1/2$, as we proved in Section 7 of \cite{LanguascoZ2015a}.
Another positive indication of this fact is Theorem 2 of Goldston-Yang \cite{GoldstonY2016},
see \eqref{Goldston-Yang} before, since assuming RH, 
they were able to reach $k=1$ (even if they do not 
give an explicit expression for  such a term, but just an estimate).

\section{Corrigendum of \cite{LanguascoZ2013a}
 and of \cite{LanguascoZ2015a}}
 \label{corrigendum}

 We take this occasion to correct the oversight occurred in Lemma \ref{Linnik-lemma2}
 which affects both  \cite{LanguascoZ2013a}
 and  \cite{LanguascoZ2015a}. The amended final statements are written before
 as Theorems  \ref{CesaroG-thm}-\ref{CesaroHL-thm}.  
 We insert here both the new parts using the original sections and equations numbering.
 
\paragraph{\textbf{Corrigendum of  \cite{LanguascoZ2013a}}}
 
In \cite{LanguascoZ2013a} the affected parts are
 Lemma 1 which has to be replaced by Lemma
\ref{Linnik-lemma2} as we said before, a part of   page 571 and  
the first part of Section 9 there. These parts 
should be substituted with what follows.

\setcounter{section}{8}
\setcounter{equation}{36} 


\smallskip
On page 571 the affected part is the estimate of
\[
\int_{(a)}
    \vert E(a,y)\vert  \, \vert e^{N z}\vert  \, \vert z\vert ^{- k - 1} \vert \omega_2(z)\vert  \, \vert \dx z\vert.
\]
Taking $a=1/N$, using   (12)-(13)
and Lemma \ref{Linnik-lemma2} before, this integral is
\begin{align*}
&\ll
 a^{- 1 / 2} e^{Na} 
 \Bigl( 
 \int_{-a}^{a} a^{{-k-1}} \dx y 
 +
  \int_{a}^{+\infty} y^{{-k-1}} \dx y 
 +
 \int_{a}^{+\infty} y^{{-k-1/2}} \log^2(y/a)\,  \dx y 
 \Bigr)
 \\
&\ll_k
  e^{N a} a^{-1/2}
  \Bigl(  a^{-k} + a^{-k+1/2}
 \int_{1}^{+\infty}v^{{-k-1/2}} \log^2 v \, \dx v
   \Bigr)
\ll_k N^{k+1/2},
\end{align*}
provided that $k>1/2$.
The previous estimate  has to be inserted in (15) and (26) and causes 
the change in the final error term of 
Theorem \ref{CesaroHL-thm}.

Moreover the first part of Section 9 should be replaced by the following.
We need $k>1$ in this section.
We first have to  establish the convergence of
\begin{equation}
\label{conv-integral-J5}
  \sum_{\ell \ge 1}
    \int_{(1/N)} 
         \vert \sum_{\rho} \Gamma(\rho)z^{- \rho} \vert
    \vert e^{N z}\vert  \vert z\vert ^{- k - 3/2}
    e^{-\pi^2 \ell^2 \Re(1/z)}  \, \vert \dx z\vert .
\end{equation}
Using the Prime Number Theorem and 
(29),
we first remark that
\begin{align}
\notag
  \vert \sum_{\rho} z^{-\rho} \Gamma(\rho) \vert 
  &=
  \vert  \Stilde(z) - \frac 1z -E(y,\frac{1}{N}) 
  \vert 
  \ll
  N   
  +
  \vert   E(y,\frac{1}{N}) \vert
\\
\label{sum-over-rho-new-HL}
  &\ll
  N + \vert z\vert ^{1/2}
    \log^2 (2N\vert y\vert) .
\end{align}

By (33)
and \eqref{sum-over-rho-new-HL}, we can write that
the  quantity in \eqref{conv-integral-J5} is
\begin{align}
\notag 
\ll N
  \sum_{\ell \ge 1}
    \int_{0}^{1/N}  
    \frac{e^{- \ell^2 N}}{\vert z\vert ^{k + 3/2}}  \,&  \dx y
    +
 N\sum_{\ell \ge 1}
    \int_{1/N}^{+\infty}   
    \frac{e^{- \ell^2 / (Ny^2)} }{\vert z\vert ^{k + 3/2}}  \, \dx y
    \\ 
    \label{HL-split}
    &
    +
    \sum_{\ell \ge 1}
    \int_{1/N}^{+\infty}   \log^{2}(2Ny)
    \frac{e^{- \ell^2 / (Ny^2)}}{\vert z\vert ^{k + 1}}   \, \dx y
    =
V_1+V_2+V_{3},
\end{align}
say.
$V_{1}$  can be estimated exactly as $U_{1}$
in Section 8
and we get $V_{1} \ll_{k} N^{k+1}$.
For $V_{2}$ we can work analogously to $U_{2}$ thus obtaining
\begin{align}
\notag
V_2&
\ll
 N \sum_{\ell \ge 1}
    \int_{1/N}^{+\infty} 
    \frac{e^{- \ell^2 / (Ny^2)}}{y^{k + 3/2}}  \, \dx y
  \ll
N^{k/2+5/4} \sum_{\ell \ge 1}
\frac{1}{\ell^{k+1/2}}
    \int_{0}^{\ell^2 N} u^{ k/2 - 3/4} 
    e^{- u}  \, \dx u 
\\
\notag
&\leq
\Gamma \Bigl( \frac{2k+1}{4} \Bigr)
N^{k/2+5/4}  \sum_{\ell \ge 1}
\frac{1}{\ell^{k+1/2}}
\ll_{k}
N^{k/2+5/4},
\end{align}
provided that $k>1/2$,
where  we used the substitution
$u=\ell^2 / (Ny^2)$.
Hence, we have 
\begin{equation}
\label{V1-V2-estim}
V_{1} + V_{2} \ll_{k} N^{k+1},
\end{equation}
provided that  $k>  1/2$.
Using the substitution $u=\ell^2 / (Ny^2)$,
we obtain
\begin{align}
\notag
V_3
&\ll
 \sum_{\ell \ge 1}
    \int_{1/N}^{+\infty} \log^{2}(2Ny) 
    \frac{e^{- \ell^2 / (Ny^2)}}{y^{k + 1}}   \, \dx y
\\
\notag
&=
\frac{N^{k/2}}{8} \sum_{\ell \ge 1}
\frac{1}{\ell^k}
    \int_{0}^{\ell^2 N} u^{ k/2 - 1} 
    \log^{2}\Bigl( \frac{4\ell^{2}N}{u} \Bigr)
    e^{- u}  \, \dx u.
\end{align}
Hence, a direct computation shows that
\begin{align}
\notag
V_3
&\ll
N^{k/2} \sum_{\ell \ge 1}
\frac{\log^{2}(\ell N)}{\ell^k}
    \int_{0}^{\ell^2 N} u^{ k/2 - 1}     e^{- u}  \, \dx u
 \\
 \notag
 & \hskip2cm
    +
N^{k/2} \sum_{\ell \ge 1}
\frac{1}{\ell^k}
     \int_{0}^{\ell^2 N} u^{ k/2 - 1}  \log^{2} (u)\,
    e^{- u}  \, \dx u
\\
\label{V3-estim}
&\ll_{k}
\Gamma(k/2)
N^{k/2} \sum_{\ell \ge 1}
\frac{\log^{2}(\ell N)}{\ell^k}
    +
N^{k/2}
\ll_{k}
N^{k/2}\log^{2} N
\end{align}
provided that $k>1$.
Inserting \eqref{V1-V2-estim}-\eqref{V3-estim}
into \eqref{HL-split} we get, for $k>1$, that
the quantity in \eqref{conv-integral-J5}
is $\ll N^{k+1}$.

The remaining part of Section 9 of \cite{LanguascoZ2013a} is left untouched.
 
\paragraph{\textbf{Corrigendum of  \cite{LanguascoZ2015a}}}
In \cite{LanguascoZ2015a} 
the affected parts are Lemma  4.1 which has to be replaced by Lemma
\ref{Linnik-lemma2} as we said before, a part of  page 1949 and 
the first part of Section 6 there. These parts 
should be substituted with what follows.
 
\setcounter{section}{6}
\setcounter{equation}{0}
\numberwithin{equation}{section}
\smallskip
On page 1949 the affected part is the following.
Recalling (2.4) and Lemma \ref{Linnik-lemma2}, we have for $a = 1 / N$ that
\begin{align*}
  \int_{(a)} \bigl \vert E(y,a) \bigr\vert^2
    \vert e^{N z} \vert \, \vert z \vert^{- k - 1} \, \vert \dx z \vert  
  &\ll_k
  e^{N a} a^{- k} 
  \ll_k N^{k} 
\end{align*}
 for $k > 1$.
For $a = 1 / N$, by (2.4) and Lemma \ref{Linnik-lemma2},
the second remainder term in (3.1) gives a contribution 
\begin{align*}
  &\ll
  N
  \int_{(1/N)} \vert E(y,1/N) \vert\vert e^{N z}\vert
    \vert z\vert ^{- k - 1} \, \vert \dx z \vert  
  \ll_k
  N^{k + 1}.
\end{align*}
This  term has to be inserted in (3.2)-(3.3) and causes the change in the error term
of Theorem \ref{CesaroHL-thm}.


\smallskip
Moreover the first part of Section 6 should be replaced by the following.
We need $k>1$ in this section.
Arguing as in Section 5
we first need to establish the
convergence of
\begin{equation}
\label{conv-integral-2}
  \sum_{\rho_{1}}
    \vert \Gamma(\rho_{1})\vert  
    \int_{(1/N)}  
      \vert \sum_{\rho_{2}} \Gamma(\rho_{2})z^{- \rho_{2}} \vert
      \vert e^{N z} \vert \,
      \vert z \vert ^{- k - 1} \, \vert z^{- \rho_{1}} \vert \, \vert \dx z \vert.
\end{equation}

Using the Prime Number Theorem and Lemma \ref{Linnik-lemma2}
 we first remark
that
\begin{align} 
  \Bigl\vert \sum_{\rho} z^{-\rho} \Gamma(\rho) \Bigr\vert  
\label{sum-over-rho-new}
  &\ll
  N + \vert z\vert ^{1/2}
    \log^2 (2N\vert y\vert).
  \end{align}
By symmetry, we may assume that $\gamma_{1}> 0$.
By \eqref{sum-over-rho-new}, 
(2.4) and (4.1),
for $y \in(-\infty, 0]$  we are first led to estimate
\begin{align*}
  \sum_{\rho_{1} \colon \gamma_{1} > 0} 
    \gamma_{1}^{\beta_{1} - 1/2}
    \exp\Bigl( -\frac \pi2  \gamma_{1}  \Bigr)
  &\Bigl(
    \int_{-1/N}^{0}  N^{ k +  2 + \beta_{1}} \, \dx y
    +
    N\int_{-\infty}^{-1/N}  
     \frac{\dx y}{\vert y\vert^{k + 1+\beta_{1}}}    
    \\
  &+
  \int_{-\infty}^{-1/N}  \log^{2} (2N \vert y \vert)
    \frac{\dx y}{\vert y\vert^{k + 1/2+\beta_{1}}}     
  \Bigr)
  \ll_k
  N^{k + 2},
\end{align*}
by the same argument used in the proof of Lemma 4.3
with $\alpha=k+1/2$ and $a=1/N$.
On the other hand, for $y > 0$ we split the range of integration into
$(0, 1 / N] \cup (1 / N, +\infty)$.
By \eqref{sum-over-rho-new}, (2.4) and Lemma 4.3
with $\alpha=k+1$ and $a=1/N$, on  $[0,1/N]$ we have
\begin{align*}
N  \sum_{\rho_{1} \colon \gamma_{1} > 0}
   \gamma_{1} ^{\beta_{1} - 1/2}
    \int_0^{1 / N} \exp\Bigl( \gamma_{1} (\arctan(N y) -\frac \pi2) \Bigr)
    \frac{\dx y}{\vert z\vert^{k + 1+\beta_{1}}}   
  &\ll_k
  N^{k + 2}.
\end{align*}
On the other interval, again by (2.4),
we have to estimate
\begin{align*}
  &
  \sum_{\rho_{1} \colon \gamma_{1} > 0}
   \gamma_{1} ^{\beta_{1} - 1/2}
    \int_{1 / N}^{+\infty} \exp\Bigl( - \gamma_{1} \arctan\frac 1{N y} \Bigr) 
      \frac{ N+ y^{1/2}\log^{2} (2Ny)}{y^{k + 1+\beta_{1}}} \dx y
 \\
  &=
  N^{k}
  \sum_{\rho_{1} \colon \gamma_{1} > 0}� \!\!\!\!
  N^{\beta_{1} }
   \gamma_{1} ^{\beta_{1} - 1/2}
    \int_{1 }^{+\infty} \exp\Bigl( - \gamma_{1} \arctan\frac 1{u} \Bigr)
      \frac{N + u^{1/2}N^{-1/2}\log^{2} (2u)}{u^{k + 1+\beta_{1}}} \dx u.
\end{align*}
Lemma 4.2
with $\alpha=k+1/2$ shows that the last
term is $\ll_k N^{k + 2}$.
This implies that the integral in 
\eqref{conv-integral-2} 
is $\ll_k N^{k + 2}$ 
provided that $k > 1$ and hence we can exchange the first summation 
with the integral in this case.

The remaining part of Section 6 of \cite{LanguascoZ2015a} is left untouched.

\vspace{0.5cm} \indent {\it
A\,c\,k\,n\,o\,w\,l\,e\,d\,g\,m\,e\,n\,t\,s.\;}  I wish to thank  the anonymous
referee for his comments and suggestions. I also wish  to thank
Alessandro Zaccagnini  and the organisers of the ``Terzo Incontro Italiano
di Teoria dei Numeri''.

\bigskip
\begin{center}

\end{center}

\bigskip
\bigskip
\begin{minipage}[t]{10cm}
\begin{flushleft}
\small{
{\sc Alessandro Languasco}
\\Universit\`a di Padova,
\\Dipartimento di Matematica
\\ Via Trieste, 63
\\ Padova, 35121, Italy
\\e-mail: languasc@math.unipd.it
}
\end{flushleft}
\end{minipage}

\end{document}